# Nontrivial Galois Module Structure of Cyclotomic Fields

## To Appear: Mathematics of Computation

Marc Conrad

Faculty of Technology, Southampton Institute
East Park Terrace, Southampton, S014 0YN Great Britain

and

Daniel R Replogle

Department of Mathematics and Computer Science
College of Saint Elizabeth, 2 Convent Road, Morristown, NJ 07960

Abstract: We say a tame Galois field extension $L/K$ with Galois group $G$ has trivial Galois module structure if the rings of integers have the property that $\mathcal{O}_L$ is a free $\mathcal{O}_K[G]$-module. The work of Greither, Replogle, Rubin, and Srivastav shows that for each algebraic number field other than the rational numbers there will exist infinitely many primes $l$ so that for each there is a tame Galois field extension of degree $l$ so that $L/K$ has nontrivial Galois module structure. However, the proof does not directly yield specific primes $l$ for a given algebraic number field $K$. For $K$ any cyclotomic field we find an explicit $l$ so that there is a tame degree $l$ extension $L/K$ with nontrivial Galois module structure.

Correspondence Address:
   Daniel R Replogle
   College of Saint Elizabeth
   Department of Mathematics and Computer Science
   2 Convent Road
   Morristown, NJ 07960

Section 1: Introduction to Cyclotomic
Swan Subgroups and Galois Module Theory

Let $G$ be a group of finite order $m$. Let $L/K$ be a tame (i.e. at most tamely ramified) Galois extension of algebraic number fields with finite Galois group $Gal(L/K) \cong G$. Let $\mathcal{O}_L$ and $\mathcal{O}_K$ denote the respective rings of algebraic integers. We say $L/K$ has a trivial Galois module structure if $\mathcal{O}_L$ is a free $\mathcal{O}_K[G]$-module. Equivalently, one says in this case $L/K$ has a normal integral basis. (Note: One always has that $\mathcal{O}_L$ is a rank one locally free $\mathcal{O}_K[G]$-module whenever $L/K$ is a tame Galois extension with Galois group $G$).

The classical Hilbert-Speiser theorem proves that any abelian extension of $\mathbb{Q}$, the field of rational numbers, has a trivial Galois module structure. Call a field $K$ Hilbert-Speiser if each tame abelian extension has a trivial Galois module structure. In [2] tame elementary abelian extensions and Swan modules are considered to find conditions a Hilbert-Speiser field must satisfy. Let $V_l = (\mathcal{O}_K/l\mathcal{O}_K)^*/Im(\mathcal{O}_K^*)$ where for any ring $S$ we let $S^*$ denote its group of multiplicative units and $Im$ denotes the image of $\mathcal{O}_K^*$ under the canonical surjection $\psi: \mathcal{O}_K \longrightarrow \mathcal{O}_K/l\mathcal{O}_K$. Then we have the following theorem.

**Theorem 1.0 ([4, Theorem 1]).** *Let $K$ be a Hilbert-Speiser number field. Then:*
  *(i) The class number of $K$ is one.*
  *(ii) For each odd prime $l$ the group $V_l$ has exponent dividing $\frac{(l-1)^2}{2}$.*
  *(iii) The group $V_2$ is nontrivial.*

This theorem and a Galois theoretic argument are used to show for any $K \neq \mathbb{Q}$ there is some odd prime $l$ for which condition (ii) is violated. Thus we have the following theorem.

**Theorem 1.1 ([4, Theorem 2]).** *Among all algebraic number fields only the rational numbers are Hilbert-Speiser.*

For the convenience of the reader we outline the ideas in the proofs of these two results. The over all idea is that $V_l^{l-1}$ will be seen to be a lower bound on the Swan subgroup, and $V_l^{(l-1)^2/2}$ will be seen to be a lower bound on the group of realisable classes. Let $\Lambda = \mathcal{O}_K[G]$. Then $\Lambda$ is an order in the group algebra $K[G]$. For each $s$ in $\mathcal{O}_K$ so that $s$ and $m$ are relatively prime one defines the Swan module $\langle s, \Sigma \rangle$ by $\langle s, \Sigma \rangle = s\Lambda + \Lambda\Sigma$, where $\Sigma = \sum_{g \in G} g$. It is easily shown each $\langle s, \Sigma \rangle$ is a rank one locally free $\Lambda$-module. (See [15, Proposition 2] for example). Hence each Swan module determines a class in the locally free classgroup $Cl(\Lambda)$. We denote the class of $\langle s, \Sigma \rangle$ by $[s, \Sigma]$. Denote the set (at this point) of all classes of Swan modules over $\Lambda$ by $T(\Lambda)$. Let $D(\Lambda)$, the kernel group, denote the subgroup of $Cl(\Lambda)$ of all classes that become trivial upon extension of scalars to the maximal order in $K[G]$ containing $\Lambda$. Let $\Gamma = \Lambda/\Sigma\Lambda$, $\psi: \mathcal{O}_K \longrightarrow \mathcal{O}_K/m\mathcal{O}_K = \overline{\mathcal{O}_K}$ denote the canonical quotient map, $\epsilon: \Lambda \longrightarrow \mathcal{O}_K$ denote the augmentation map, and let $\bar{\epsilon}: \Gamma \longrightarrow \overline{\mathcal{O}_K}$ be induced from $\epsilon$. The main result of [9] shows that if the group algebra $K[G]$ satisfies an Eichler condition (which is always satisfied if $G$ is abelian) then there exists an exact Mayer-Vietoris sequence:

$$\mathcal{O}_K^* \times \Gamma^* \xrightarrow{h} \overline{\mathcal{O}_K}^* \xrightarrow{\delta} D(\Lambda) \longrightarrow D(\Gamma) \oplus D(\mathcal{O}_K) \longrightarrow 0, \tag{1.0}$$

where for any ring $S$ we have been denoting its group of multiplicative units by $S^*$. It is further shown in [9] that the map $h$ is given by: $h[(u,v)] = \psi(u)\overline{\epsilon}(v)^{-1}$. From [15] we have that the map $\delta$ is given by: $\delta(s) = [s, \Sigma]$. Hence we have from (1.0) that $T(\Lambda)$ is a subgroup of $D(\Lambda)$ and :

$$T(\Lambda) \cong (\overline{\mathcal{O}_K})^*/h(\mathcal{O}_K^* \times \Gamma^*). \tag{1.1}$$

Now let $V_l = (\mathcal{O}_K/l\mathcal{O}_K)^*/Im(\mathcal{O}_K^*)$ and set $G \cong C_l$, the cylcic group of prime order $l$. Then in [4] a lower bound on $T(\Lambda)$ is deduced from (1.1) above. Namely one may show for $G \cong C_l$ and $\Lambda = \mathcal{O}_K[C_l]$ there is a natural surjective map:

$$T(\Lambda) \longrightarrow V_l^{l-1}. \tag{1.2}$$

See [4,Theorem 5] which handles the elementary abelian group case.

Now consider tame Galois extensions of $K$ with $Gal(L/K) \cong C_l$ where $C_l$ is the cyclic group of order $l$. Of course each such extension $L$ has the property that $\mathcal{O}_L$ is a locally free rank one $\mathcal{O}_K[G]$-module. Hence each $\mathcal{O}_L$ determines a tame Galois module class $[\mathcal{O}_L]$ in the locally free classgroup $Cl(\Lambda)$. Denote by $R(\Lambda)$ the set of all such classes. Observe $R(\Lambda)$ measures to what extent a given base field $K$ fails to have Galois extensions with Galois group $G$ having nontrivial Galois module structure. In [7] McCulloh gives an explicit description of $R(\Lambda)$ as a subgroup of $Cl(\Lambda)$ for the case considered here. We note in [8] McCulloh shows that $R(\Lambda)$ is in fact a subgroup of $Cl(\Lambda)$ for all abelian $G$, however that result is not used here. Using the description of $R(\Lambda)$ from [7] we have the following relationship between $R(\Lambda)$, $D(\Lambda)$, and $T(\Lambda)$ stated for the case we are considering. If $G \cong C_l$ where $l > 2$ is prime, then by [4, Proposition 4] the Swan subgroup has the property that

$$T(\Lambda)^{(l-1)/2} \subseteq R(\Lambda) \cap D(\Lambda). \tag{1.3}$$

Observe that condition (ii) of Theorem 1.0 follows from (1.2) and (1.3). The proof of Theorem 1.1 uses the Cebotarev density theorem to establish the existence of infinitely many primes for which condition (ii) of Theorem 1.0 is violated. It is not clear whether one can use the proof of Theorem 1.1 to actually write down any explicit such primes. In fact as the proof of [4, Theorem 2] uses the Cebotarev density theorem to show the existence of such a prime by showing infinitely many must exist it is unlikely one could use that argument to write down any such primes. We adopt the following terminology. If given a field $K$ we have an explicit prime $l$ so that there is a tame Galois extension $L/K$ with $Gal(L/K) \cong C_l$, the cyclic group of order $l$, so that $\mathcal{O}_L$ is not a free $\mathcal{O}_K[C_l]$-module we say $K$ is not Hilbert-Speiser for $l$. Then Theorem 1.1 says for any algebraic number field $K$ there is some $l$ for which $K$ is not Hilbert-Speiser for $l$.

In this article we show that for any cyclotomic field one may in fact find a specific prime $l$ for which the field is not Hilbert-Speiser for $l$. That is let $K_n = \mathbb{Q}(\zeta_n)$ where $\zeta_n$ is a primitive $n$th root of unity. Note without loss of generality we may assume $n \not\equiv 2 \mod 4$ as if $n \equiv 2 \mod 4$ then $K_n = K_{\frac{n}{2}}$.

For $K_n$ the ring of algebraic integers is $\mathbb{Z}[\zeta_n]$ which we denote by $\mathcal{O}_n$. The class number of an algebraic number field $K$ is denoted by $h_K$ and for $K_n$ is denoted by $h_n$. In light of Theorem 1.0 to accomplish our goal it suffices to prove the following theorem, the main result for this article.

**Theorem 1.2 (Main Theorem).** *Let $K_n$ be as above and assume $h_n = 1$. One may find an explicit prime $l$ so that $V_{l,n} = (\mathcal{O}_n/l\mathcal{O}_n)^*/Im(\mathcal{O}_n^*)$ does not have exponent dividing $\frac{(l-1)^2}{2}$. Hence for every cyclotomic field one may find an explicit prime $l$ so that there is a tame degree $l$ field extension with nontrivial Galois module structure.*

We begin our study by taking care of two special cases. We first treat the case $h_n \neq 1$ by providing some details about condition (i) of Theorem 1.0. Using [6] one can show any field of class number not equal to one has a quadratic extension which does not posses a relative integral basis. Thus we have the following restating this in our language.

**Theorem 1.3.** *Let $K$ be an algebraic number field with class number $h_K \neq 1$. Then $K$ is not Hilbert-Speiser for $l = 2$.*

While this result is decisive for fields of class number not equal to 1, we must note much more recent studies of Galois Module Theory have been made. For example the work of Fröhlich, [3], gives an excellent account of Galois module classes. The results in [1] which we state below provide an important Galois module structure result for extensions of cyclotomic fields.

Next we treat the case when $n$ is prime. We note for this case one need not make any assumptions on the class number.

**Propostion 1.4.** *For $n > 3$ prime one has $K_n$ is not Hilbert-Speiser for $n$.*

*Proof.* Let $n > 3$ be prime. We show $V_n = (\mathcal{O}_n/n\mathcal{O}_n)^*/Im(\mathcal{O}_n)$ is an elementary abelian group of exponent $n$, from which the result follows from Theorem 1.0 (ii). As $n$ totally ramifies, we have $(\mathcal{O}_n/l\mathcal{O}_n)^* \cong C_{n-1} \times C_n^{n-2}$, and from Dirichlet's unit theorem we have $\mathcal{O}_n^* \cong \langle -\zeta_n \rangle \times \langle \epsilon_1 \rangle \times \cdots \times \langle \epsilon_{(n-3)/2} \rangle$, where the $\epsilon_i$ are a system of fundamental units. Since $\zeta_n$ is not congruent to 1 mod $l$, $V_n^{n-1} \cong C_n^{l-2-j}$ for some integer $j$ with $1 \leq j \leq (n-1)/2$. Now as $n \geq 5$, we have $V_n^{n-1}$ is a nontrivial elemetrary abelian group of exponent $n$. □

Notes: 1. Proposition 1.4 is essentially [10, Proposition 15]. 2. In the case both $n$ is prime and $n \nmid h_n^+$, where $h_n^+$ is the class number of the maximal real subfield of $K_n$ one may in fact explicitly compute $T(\mathbb{Z}[\zeta_n]C_n)$. See [11, Theorem 1]. 3. For applications to Hopf orders in $\mathbb{Z}[\zeta_n]C_n$ when $n$ is a prime power, see [12].

So, we need to only consider the case where either $n = 3$ or $n$ is composite so that $h_n = 1$. This reduces the question to the following list. (See [16,Theorem 11.1], for example, for a listing of those cyclotomic fields $K_n$ with $n \not\equiv 2 \mod 4$ so that $h_n = 1$.)

**List 1.5.** *Let n not be congruent to 2 mod 4 and not a prime greater than 3. Then $h_n = 1$ only if $n = 3, 4, 8, 9, 12, 15, 16, 20, 21, 24, 25, 27, 28, 32, 33, 35, 36, 40, 44, 45, 48, 60,$ or $84$.*

In the next section we will provide several primes $l$ for each $n$ in this list so that $K_n$ is not Hilbert-Speiser for $l$. That is for each $n$ we will exhibit several primes $l$ so that $V_{l,n}$ does not have exponent dividing $\frac{(l-1)^2}{2}$. This will complete the proof of the theorem.

Before beginning we introduce the important related result of Chan and Lim, [1]. For any Galois extension of number fields $L/K$ with Galois group $G$ one may define the associated order by $\mathcal{A}_{L/K} = \{\alpha \in K[G] | \alpha \mathcal{O}_L \subseteq \mathcal{O}_L\}$. Then a classical result is $L/K$ is at most tamely ramified if and only if $\mathcal{A}_{L/K} = \mathcal{O}_K[G]$. So, one may generalize the notion of trivial Galois module structure to say a Galois extension of number fields $L/K$ has trivial Galois module structure if $\mathcal{O}_L$ is a free $\mathcal{A}_{L/K}$-module. Chan and Lim show, [1], any Galois extension $L/K$ with $L$ and $K$ both cyclotomic fields has a trivial Galois module structure. Hence, for the primes $l$ for which we detect the existence of degree $l$ extensions with nontrivial Galois module structure we know these exension fields are not cyclotomic fields.

## Section 2: Cyclotomic Units and Nontrivial Galois Module Structure

To complete the proof of Theorem 1.2 it suffices to exhibit a prime $l$ for each $n$ in List 1.5 so that $V_{l,n} \cong (\mathcal{O}_n/l\mathcal{O}_n)^*/Im(\mathcal{O}_n^*)$ is not of exponent dividing $\frac{(l-1)^2}{2}$. Of course the issues are finding an $l$ that splits nicely in $\mathcal{O}_n$; then representing the quotient $\mathcal{O}_n/l\mathcal{O}_n$ in a way amenable to computation; then finding the image of the units of $\mathcal{O}_n$ in this representation; then last developing a computer algorithm to solve the problem. We will not use the results in [5] here, however we must note our method is somewhat analogous.

The first two parts of our task are completed by the following lemma. Let $\mathbb{F}_l$ denote the field of $l$ elements and $\mathbb{F}_l[x]$ its polynomial ring. Further let $\phi$ be the usual Euler $\phi$-function.

**Lemma 2.0.** *For an $l$ that does not divide $n$ let $f$ be minimal such that $n | k := l^f - 1$, and $g = \phi(n)/f$. Then $(\mathcal{O}_n/l\mathcal{O}_n)^* \cong (\mathbb{Z}/k\mathbb{Z})^g$, where the isomorphism is explicitly given by the following maps $\kappa$, $\lambda_i$, and $\mu_i$:*

$$\kappa : \mathcal{O}_n/l\mathcal{O}_n \longrightarrow \prod_{i=1}^{g} \mathcal{O}_n/P_i\mathcal{O}_n \quad (2.0)$$

$$a + l\mathcal{O}_n \mapsto (a + P_i\mathcal{O}_n)_{i=1,\ldots,g},$$

*where each $P_i$ is an ideal generated by two generators $P_i = \langle l, p_i \rangle$, and $p_i$ is a polynomial of degree $f$, which we obtain by factoring the n-th cyclotomic polynomial $\Phi_n(x)$ as a polynomial in $\mathbb{F}_l[x]$ into irreducible polynomials.*

$$\lambda_i : \mathcal{O}_n/P_i\mathcal{O}_n \to \mathbb{F}_l[x]/(p_i) \quad (2.1)$$

$$a + P_i\mathcal{O}_n \mapsto a + (p_i)$$

for $i = 1, \ldots, g$.

$$\mu_i : (\mathbb{F}_l[x]/(p_i))^* \to \mathbb{Z}/k\mathbb{Z} \qquad (2.2)$$

$$a + (p_i) \mapsto b \bmod k.$$

where $a + (p_i) = w_i^b$ and $w_i$ is a (fixed) generator of the cyclic group $(\mathbb{F}_l[x]/(p_i))^*$.

*Proof.* The isomorphism (2.0) follows by [16, Theorem 2.13], which describes the splitting behaviour of a prime in $K_n$ and the Chinese remainder theorem for ideals in an algebraic number field. The isomorphism (2.1) can be proved straightforwardly as $\mathcal{O}_n \cong \mathbb{Z}[x]/(\Phi_n(x))$ and (2.2) is the well known fact that the multiplicative group of a finite field is cyclic. □

The larger issue is the units. We show first that it is sufficient to consider cyclotomic units. We introduce some mostly standard notation at this point: Let $W_n$ denote the roots of unity in $\mathcal{O}_n$. Let $E_n$ denote the units of $\mathcal{O}_n$. Last let $E_n^+$ denote the units in the ring of algebraic integers of the maximal real subfield of $K_n$. Then we have $[E_n : W_n E_n^+] = Q$ where $Q$ equals one if $n$ is a prime power and equals two otherwise. For details see chapters 4 and 8 of [16].

Let $C_n^+$ denote the cyclotomic units in $E_n^+$. Let $h_n^+$ denote the class number of the maximal real subfield of $K_n$. Next we have the following result due to Sinnott [14] stated only for what we are considering here.

**Theorem 2.1 [14].** *The group $E_n^+/C_n^+$ is finite and $[E_n^+ : C_n^+] = 2^b h_n^+$. The integer $b$ is defined by $b = 0$ if $r = 1$ and $b = 2^{r-2} + 1 - r$ if $r > 1$, where $r$ is the number of distinct primes dividing $n$.*

This theorem has two immediate corollaries which we prove.

**Corollary 2.2.** *If $h_n = 1$ then $E_n^+ = C_n^+$.*

*Proof.* For each $K_n$ of class number one three or fewer primes divide $n$. Hence $b = 0$ in Theorem 2.1. □

**Corollary 2.3.** *For each $n$ of List 1.4 we have $C_n = E_n$.*

*Proof.* Note first that $E_n^+ = C_n^+$ by Corollary 2.2 and obviously we have $C_n \supseteq W_n C_n^+$. If $n$ is a prime power we have $Q = 1$ and therefore $E_n = W_n E_n^+ = W_n C_n^+ = C_n$.

In the case $Q = 2$ the result is proven when we show that $[C_n : W_n C_n^+] > 1$. This follows because $1 - \zeta_n \in C_n$ but $\zeta_n^\nu(1 - \zeta_n) = \zeta_n^\nu - \zeta_n^{\nu+1} \notin \mathbb{R} \supseteq C_n^+$ for each $\nu = 1, \ldots, n$. (A necessary condition for $\zeta_n^\nu - \zeta_n^\mu \in \mathbb{R}$ is $2(\nu + \mu) = n$ which is impossible for $\mu = \nu + 1$ as $n \not\equiv 2 \bmod 4$.) □

From this it follows we may consider just cyclotomic units. That is, the group of cyclotomic units generates the full group of units. From [2] we have an explicit description of a basis of the cyclotomic units in $K_n$.

**Proposition 2.4 [2].** *A basis of the group of cyclotomic units is explicitely given as a subset of the set*

$$\{\frac{1-\zeta_q^a}{1-\zeta_q} \mid q|n,\ 1 < a < q,\ (a,q) = 1,\ q \text{ is a prime power}\}$$

$$\cup \{1 - \zeta_d^a \mid d|n,\ 1 < a < d,\ (a,d) = 1,\ d \text{ is not a prime power}\} \qquad (2.3)$$

$$\cup \{\pm\zeta_n\}.$$

An algorithm which computes a basis according to [2] is implemented in SIMATH [13]. Now we work toward developing an algorithm to compute $V_{l,n}$. We will need the following lemma.

**Lemma 2.5.** *Let $k, s, g \in \mathbb{N}$ and $G = (\mathbb{Z}/k\mathbb{Z})^g$. Further let $H = \langle v_1, \ldots, v_g \rangle$ be a subgroup of $G$ which is generated by the $g$-tuples $v_i = (v_{1,i}, \ldots, v_{g,i})$ and $V = (v_{i,j})_{1 \leq i,j \leq g}$ the $g \times g$-matrix which is made up by the $v_i$. Then we have:*

*If there exists a prime $p$ such that $p|k, p|\det(V)$ and $p \nmid s$, then the exponent of $G/H$ does not divide $s$.*

*Proof.* The exponent of $G/H$ does not divide $s$ means that $\exists w \in G$, such that $sw \notin H$, or – more concretely – $\exists w \in G$ with $\forall \alpha := (\alpha_1, \ldots, \alpha_g) \in \mathbb{Z}^g$, such that

$$\sum_{i=1}^{g} \alpha_i v_i \not\equiv sw \bmod k. \qquad (2.4)$$

A sufficient condition for (2.4) is obviously that there exists a prime $p|k$ with

$$\sum_{i=1}^{g} \alpha_i v_i \not\equiv sw \bmod p. \qquad (2.5)$$

The only primes $p$ which are in question are those $p$ with $p \nmid s$, otherwise $(0, \ldots, 0)$ would be a solution. For those $p$ the number $s^{-1}$ is well defined mod $p$, and we obtain from (2.5) that we have to show: $\exists w \in G$ with $\alpha s^{-1} V \not\equiv w \bmod p$. So in fact we have to show that the endomorphism

$$\varphi: (\mathbb{Z}/p\mathbb{Z})^g \to (\mathbb{Z}/p\mathbb{Z})^g$$

$$x \mapsto xs^{-1}V$$

is not surjective, which is equivalent to $\det(s^{-1}V) \equiv 0 \bmod p$, or $p|\det(V)$. □

**Remark 2.6** Of course the subgroup $H$ can be generated by fewer then $g$ vectors. Then we obtain the lemma by dropping the condition $p|\det(V)$.

**Remark 2.7** Note that we do not need to factor $k$ in order to check whether a $p$ as required in Lemma 2.5 exists. We can use the following algorithm:
i) Let $c := \gcd(k, \det(V))$.
ii) do $t := c;\ c := c/\gcd(c, s)$ while $c \neq t$.

iii) if $c = 1$ then $\not\exists p$ else $\exists p$.

Except for the proof that the group of cyclotomic units already is the full unit group we have not used the fact that $h_n = 1$. So we define $V_{l,n}^{cyc} \cong (\mathcal{O}_n/l\mathcal{O}_n)^*/Im(C_n)$ where we have replaced the full unit group $\mathcal{O}_n^*$ by the group of cyclotomic units $C_n$. With the results of Lemmata 2.0 and 2.5 we obtain the following algorithm.

**Algorithm 2.8** If for a given pair $(l, n)$ with $l \not| n$ the following algorithm returns "yes", then the exponent of $V_{l,n}^{cyc}$ is not dividing $(l-1)^2/2$.

i) Factor the $n$-th cyclotomic polynomial $\Phi_n$ modulo $l$. Let $L = \{p_1, \ldots, p_g\}$ be the list of irreducible factors of $\Phi_n$. Let further $f$ be the degree of one of the $p_i$ (which is in fact the same for all $p_i$) and set $k = l^f - 1$.

ii) For each polynomial $p_i \in L$ determine a generator $w_i$ of the finite field $\mathbb{F}_l[x]/(p_i(x))$.

iii) Determine a (finite) set $B = \{b_1, \ldots, b_s\}$ of generators of cyclotomic units. This may be done e.g. with the function `lcyubas()` of SIMATH, [13].

iv) Compute the matrix $V \in (\mathbb{Z}/k\mathbb{Z})^{s \times g}$ where each entry $v_{i,j}$ of $V$ is given by solving the discrete logarithm problem $b_i \equiv w_j^{v_{i,j}} \mod p_j$ in $\mathbb{F}_l[x]$.

v) If $s > g$ eleminate rows of $V$ by Gauss elimination over the ring $\mathbb{Z}/k\mathbb{Z}$, until $V$ is a square matrix. If $s < g$ generate a square matrix by adding zero rows (or set $\det(V) = 0$ in the next step).

vi) Check (via Remark 2.7) if there exists a prime $p$ such that $p | \gcd(\det(V), k)$ and $p \not| ((l-1)^2/2)$. If such a $p$ exists return "yes".

This algorithm has been implemented and is available as function `iscynfHS()` in Version 4.5 of SIMATH.

For each $n$ in List 1.5 we now in fact exhibit several primes $l$ for which $V_{l,n}^{cyc}$ is not of exponent dividing $\frac{(l-1)^2}{2}$. For these $n$ one has $V_{l,n}^{cyc} = V_{l,n}$ as $h_n = 1$. So, we provide a table obtained using the above alogrithm that for each $n$ in List 1.5 we have $K_n = \mathbb{Q}(\zeta_n)$ and for the $l$ associated to $n$ we have that there is a tame Galois field extension $L/K_n$ of degree $l$ for which $l$ splits in $K_n/\mathbb{Q}$ and the extension $L/K_n$ is so that $L/K_n$ is not Hilbert-Speiser, i.e. does not have a normal integral basis. This together with Theorem 1.3 and Proposition 1.4 proves Theorem 1.2. We note the algorithm presented applies to those $n$ treated in Proposition 1.4 also. For example for $n = 7$ one obtains the list $l = 5, 11, 13, 17, 23, 37, 41, 53, 67, 79, 83, 97$. Of course, these $l$ split whereas 7 totally ramifies in $K_7$.

| Table 2.9 | |
|---|---|
| $n$ | $l$ |
| 3 | 17, 29, 41, 53, 59, 71, 83, 89 |
| 4 | 11, 19, 23, 43, 47, 59, 67, 71, 79, 83 |
| 8 | 5, 11, 13, 19, 23, 29, 37, 43, 47, 53, 59, 61, 67, 71, 79, 83 |
| 9 | 5, 7, 11, 13, 17, 31, 43, 53, 61, 67, 71, 79, 89, 97 |
| 12 | 11, 17, 19, 23, 29, 41, 43, 47, 53, 59, 67, 71, 79, 83, 89 |
| 15 | 7, 11, 13, 17, 19, 23, 29, 37, 41, 43, 47, 53, 59, 67, 71, 73, 79, 89 |
| 16 | 3, 5, 11, 13, 19, 23, 29, 37, 41, 43, 47, 53, 59, 61, 67, 71, 73, 79, 89 |
| 20 | 7, 11, 13, 17, 19, 23, 29, 37, 43, 47, 53, 59, 67, 71, 73, 79, 89 |
| 21 | 5, 11, 13, 17, 29, 37, 41, 67, 71, 79, 83, 97 |
| 24 | 5, 11, 13, 17, 19, 23, 29, 37, 41, 43, 47, 53, 59, 61, 67, 71, 79, 83, 89 |
| 25 | 7, 11, 31, 43 |
| 27 | 17, 19, 37, 53, 73 |
| 28 | 3, 5, 11, 13, 17, 37, 41, 43, 53, 71, 83, 97 |
| 32 | 3, 5, 7, 17, 23, 41, 47, 71, 73, 79 |
| 33 | 5, 23, 43, 89 |
| 35 | 3, 11, 13, 29, 41, 43 |
| 36 | 5, 7, 11, 13, 17, 19, 53, 61, 71, 89, 97 |
| 40 | 3, 7, 11, 13, 17, 19, 23, 29, 37, 43, 47, 53, 59, 61, 67, 71, 73, 79, 89 |
| 44 | 3, 5, 23, 43, 67 |
| 45 | 11, 17, 19, 31, 37, 53, 61, 71, 73, 89 |
| 48 | 5, 11, 13, 17, 19, 23, 29, 37, 41, 43, 47, 53, 59, 71, 73, 79, 89 |
| 60 | 7, 11, 13, 17, 19, 23, 29, 37, 41, 43, 47, 53, 59, 67, 71, 73, 79, 89 |
| 84 | 5, 11, 13, 17, 29, 37, 41, 43, 71, 83, 97 |

## Acknowledgements

The authors would wish to thank L. Washington for conversations regarding this project and for his excellent text ([16]) that both use frequently.